\documentclass[]{amsart}

\usepackage{geometry}
\newtheorem{theorem}{Theorem}[section]

\newtheorem{example}[theorem]{Example}

\newtheorem{remark}[theorem]{Remark}

\newcommand{\rank}{\text{rank}}

\begin{document}

\title[Local classification and examples]{Local classification and examples of an important class of paracontact metric manifolds}

\author[V. Mart\'{\i}n-Molina]{Ver\'{o}nica Mart\'{\i}n-Molina}
\email{vmartin@unizar.es}
\address{Centro Universitario de la Defensa de Zaragoza, Academia General Militar, Ctra. de Huesca s/n, 50090 Zaragoza, SPAIN\\
 and I.U.M.A, Universidad de Zaragoza, SPAIN}

\thanks{The author is is partially supported by the PAI group
FQM-327 (Junta de Andaluc\'ia, Spain), the group Geometr\'ia E15 (Gobierno de Arag\'on, Spain),  the MINECO grant MTM2011-22621 and the ``Centro Universitario de la Defensa de Zaragoza'' grant ID2013-15.}

\begin{abstract}
We study paracontact metric $(\kappa,\mu)$-spaces with $\kappa=-1$, equivalent to $h^2=0$ but not $h=0$. In particular, we will give an alternative proof of Theorem~3.2 of \cite{mio} and present examples of paracontact metric $(-1,2)$-spaces and $(-1,0)$-spaces of arbitrary dimension with tensor $h$ of every possible constant rank. We will also show explicit examples of paracontact metric $(-1, \mu)$-spaces with tensor $h$ of non-constant rank, which were not known to exist until now.
\end{abstract}



\maketitle

\section{Introduction}

Paracontact metric manifolds, the odd-dimensional analogue of paraHermitian manifolds, were first introduced in \cite{kaneyuki} and they have been the object of intense study recently, particularly since the publication of \cite{zamkovoy}. An important  class among paracontact metric manifolds is that of the $(\kappa,\mu)$-spaces, which satisfy the nullity condition \cite{CKM}
\begin{equation}\label{kappamu}
R(X,Y)\xi=\kappa(\eta(Y)X-\eta(X)Y)+\mu(\eta(Y)hX-\eta(X) hY),
\end{equation}
for all $X,Y$ vector fields on $M$, where $\kappa$ and $\mu$ are constants and $h=\frac12 L_\xi \varphi$.

This class includes the paraSasakian manifolds \cite{kaneyuki,zamkovoy}, the paracontact metric manifolds satisfying $R(X,Y)\xi=0$ for all $X,Y$ \cite{zamkovoy-arxiv2}, certain $g$-natural paracontact metric structures constructed on unit tangent sphere bundles \cite{mio-calvaruso}, etc.

The definition of a paracontact metric $(\kappa,\mu)$-space was motivated by the relationship between contact metric and paracontact geometry. More precisely, it was proved in \cite{mino-pacific} that any non-Sasakian contact metric $(\kappa,\mu)$-space accepts two paracontact metric $(\widetilde\kappa,\widetilde\mu)$-structures with the same contact form. On the other hand, under certain natural conditions, every non-paraSasakian paracontact $(\widetilde{\kappa},\widetilde{\mu})$-space admits a contact metric $(\kappa,\mu)$-structure compatible with the same contact form (\cite{CKM}).

Paracontact metric $(\kappa,\mu)$-spaces satisfy that $h^2=(\kappa+1)\phi^2$ but this condition does not give any type of restriction over the value of $\kappa$, unlike in contact metric geometry, because the metric of a paracontact metric manifold is not positive definite. However, it is useful to distinguish the cases $\kappa>-1$, $\kappa<-1$ and $\kappa=-1$. In the first two, equation \eqref{kappamu} determines the curvature completely and either the tensor $h$ or $\varphi h$ are diagonalisable \cite{CKM}. The case $\kappa=-1$ is equivalent to $h^2=0$ but not to $h=0$. Indeed, there are examples of paracontact metric $(\kappa,\mu)$-spaces with $h^2=0$ but $h\neq0$, as was first shown in \cite{mino-kodai,CKM,CP,murathan}.

However, only some particular examples were given of this type of space and no effort had been made to understand the general behaviour of the tensor $h$ of a paracontact metric $(-1,\mu)$-space until the author published \cite{mio}, where a local classification depending on the rank of $h$ was given in Theorem~3.2. The author also provided explicit examples of all the possible constant values of the rank of $h$ when $\mu=2$. She explained why the values $\mu=0$ and $\mu=2$ are special and studying them is enough. Finally, she showed some paracontact metric $(-1,0)$-spaces of any dimension with $\rank(h)=1$ and of paracontact metric $(-1,0)$-spaces of dimension $5$ and $7$ for any possible constant rank of $h$. These were the first examples of this type with $\mu\neq2$ and dimension greater than $3$.

In the present paper, after the preliminaries section,  we will give an alternative proof of Theorem~3.2 of \cite{mio} that does not use \cite{O} and we will complete the examples of all the possible cases of constant rank of $h$ by presenting $(2n+1)$-dimensional paracontact metric $(-1,0)$-spaces with $\rank(h)=2,\ldots,n$. Lastly, we will also show the first explicit examples ever known of paracontact metric $(-1,2)$-spaces and $(-1,0)$-spaces with $h$ of non-constant rank.

\section{Preliminaries}

An \emph{almost paracontact structure} on a
$(2n+1)$-dimensional smooth manifold $M$ is given by a
$(1,1)$-tensor field $\varphi$, a vector field $\xi$ and a
$1$-form $\eta$ satisfying the following conditions \cite{kaneyuki}:
\begin{enumerate}
  \item[(i)] $\eta(\xi)=1$, \ $\varphi^2=I-\eta\otimes\xi$,
  \item[(ii)] the eigendistributions ${\mathcal D}^+$ and ${\mathcal D}^-$ of $\varphi$ corresponding to the eigenvalues $1$ and $-1$, respectively, have equal dimension $n$.
\end{enumerate}

It follows that $\varphi\xi=0$, $\eta\circ\varphi=0$ and $\rank (\varphi)=2n$.
If an almost paracontact manifold admits a semi-Riemannian metric $ g$ such that
\[
g(\varphi X,\varphi Y)=- g(X,Y)+\eta(X)\eta(Y),
\]
for all  $X,Y$ on $M$, then $(M,\varphi,\xi,\eta, g)$ is called an \emph{almost paracontact metric manifold}. Then $g$ is necessarily of signature $(n+1,n)$ and satisfies $\eta=g(\cdot,\xi)$ and ${g}(\cdot,\varphi\cdot)=-{g}(\varphi\cdot,\cdot)$.

We can now define the \emph{fundamental $2$-form} of the almost paracontact metric manifold by
$\Phi(X,Y)={g}(X,\varphi Y)$. If $d\eta=\Phi$, then $\eta$ becomes a contact form (i.e. $\eta \wedge (d\eta)^n \neq0$) and $(M,\varphi,\xi,\eta, g)$ is said to be a \emph{paracontact
metric manifold}.

We can also define on a paracontact metric manifold the tensor field ${h}:=\frac{1}{2}L_{\xi}\varphi$, which is symmetric with respect to ${g}$ (i.e. $g(hX,Y)=g(X,hY)$, for all $X,Y$), anti-commutes with $\varphi$  and satisfies $h\xi=\text{tr} h=0$ and the  identity $\nabla\xi=-\varphi+\varphi h$ (\cite{zamkovoy}).  Moreover,  it  vanishes identically if and only if $\xi$ is a Killing vector field, in which case $(M,\varphi,\xi,\eta, g)$ is called a \emph{K-paracontact manifold}.

An almost paracontact structure is said to be \emph{normal} if and only if the tensor $[\varphi,\varphi]-2d\eta\otimes\xi=0$, where $[\varphi,\varphi]$ is the Nijenhuis tensor of $\varphi$ \cite{zamkovoy}:
\[
[\varphi,\varphi](X,Y)=\varphi^2 [X,Y]+[\varphi X,\varphi Y]-\varphi [\varphi X,Y]-\varphi [X,\varphi Y].
\]
A normal paracontact metric manifold is said to be a \emph{paraSasakian manifold} and is in particular $K$-paracontact. The converse holds in dimension $3$ (\cite{calvaruso}) but not in general in higher dimensions. However, it was proved in Theorem~3.1 of \cite{mio} that it also holds for $(-1,\mu)$-spaces.
Every paraSasakian manifold satisfies
\begin{equation}\label{parasasakian-r}
R(X,Y)\xi =-(\eta(Y)X-\eta(X)Y),
\end{equation}
for every $X,Y$ on $M$. The converse is not true, since Examples~3.8--3.11 of \cite{mio} and Examples~\ref{ex-mu0-h2+} and \ref{ex-mu0-non-constant} of the present one show that there are examples of paracontact metric manifolds satisfying equation \eqref{parasasakian-r} but with $h\neq0$ (and therefore not K-paracontact or paraSasakian). Moreover, it is also clear in Example \ref{ex-mu0-non-constant} that the rank of $h$ does not need to be constant either, since $h$ can be zero at some points and non-zero in others.

The main result of \cite{mio} is the following local classification of paracontact metric $(-1,\mu)$-spaces:

\begin{theorem}[\cite{mio}]\label{th-h}
Let $M$ be a $(2n+1)$-dimensional paracontact metric $(-1,\mu)$-space. Then we have one of the following possibilities:
\begin{enumerate}
\item
either $h=0$ and $M$ is paraSasakian,

\item
or $h\neq 0$ and $\text{rank} (h_p)\in \{1,\ldots,n \}$ at every $p  \in M$ where $h_p \neq 0$. Moreover, there exists a basis $\{ \xi_p, X_1,Y_1,\ldots,X_n,Y_n \}$ of $T_p(M)$ such that
the only non-vanishing components of $g$ are
\[
g_p(\xi_p,\xi_p)=1, \quad g_p(X_i,Y_i)=\pm 1,
\]
and
\begin{equation*}
{h_p}_{| \langle X_i,Y_i \rangle}=
\begin{pmatrix}
0 & 0\\
1 & 0
\end{pmatrix}
\quad
\text{ or }
\quad
{h_p}_{| \langle X_i,Y_i \rangle}=
\begin{pmatrix}
0 & 0\\
0 & 0
\end{pmatrix},
\end{equation*}
where obviously there are exactly $\text{rank} (h_p)$ submatrices of the first type.

If $n=1$, such a basis $\{ \xi_p, X_1,Y_1\}$ also satisfies that
\[
\varphi_p X_1 =\pm  X_1, \quad \varphi_p Y_1 =\mp  Y_1,
\]
and the tensor $h$ can be written as
\begin{equation*}
{h_p}_{| \langle\xi_p, X_1,Y_1 \rangle}=
\begin{pmatrix}
0 & 0 & 0\\
0 & 0 & 0\\
0 & 1 & 0
\end{pmatrix}.
\end{equation*}
\end{enumerate}
\end{theorem}


\medskip

Many examples of paraSasakian manifolds are known. For instance, hyperboloids $\mathbb{H}^{2n+1}_{n+1}(1)$ and the hyperbolic Heisenberg group ${\mathcal H}^{2n+1}=\mathbb{R}^{2n}\times\mathbb{R}$, \cite{ivanov}. We can also obtain ($\eta$-Einstein) paraSasakian manifolds from contact $(\kappa,\mu)$-spaces with $|1-\frac{\mu}{2}|<\sqrt{1-\kappa}$. In particular, the tangent sphere bundle $T_1N$ of any space form $N(c)$ with $c<0$ admits a canonical $\eta$-Einstein paraSasakian structure,  \cite{nuestro-mino}. Finally, we can see how to construct explicitly a paraSasakian structure on a Lie group (see Example 3.4 of \cite{mio}) or on the unit tangent sphere bundle, \cite{mio-calvaruso}.

On the other hand, until \cite{mio} only some types of non-paraSasakian paracontact metric $(-1,\mu)$-spaces were known:
\begin{itemize}

\item $(2n+1)$-dimensional paracontact metric $(-1,2)$-space with $\text{rank}(h)=n$, \cite{CKM}.

\item $3$-dimensional paracontact metric $(-1,2)$-space with $\text{rank}(h)=n=1$, \cite{murathan}.

\item $3$-dimensional paracontact metric  $(-1,0)$-space with $\text{rank}(h)=n=1$. This example is not paraSasakian but it satisfies \eqref{parasasakian-r}, \cite{CP}.

\end{itemize}

The answer to why there seems to be only examples of paracontact metric $(-1,\mu)$-spaces with $\mu=2$ or $\mu=0$ is a ${\mathcal D}_{c}$-homothetic deformation, i.e. the following change of a paracontact metric structure $(M,\varphi,\xi,\eta,g)$ \cite{zamkovoy}:
\begin{equation*}
\varphi':=\varphi, \quad \xi':=\frac{1}{c}\xi, \quad \eta':=c \eta, \quad g':=c g + c(c-1)\eta\otimes\eta,
\end{equation*}
for some $c\neq0$.

It is known that $(\varphi',\xi',\eta',g')$ is again a paracontact metric structure on $M$ and that $K$-paracontact and paraSasakian structures are also preserved. However, curvature conditions like $R(X,Y)\xi=0$ are destroyed, since paracontact metric $(\kappa,\mu)$-spaces become other paracontact metric $(\kappa',\mu')$-spaces with
\begin{equation*}
\kappa'=\frac{\kappa+1-c^2}{c^2}, \quad \mu'=\frac{\mu-2+2c}{c}.
\end{equation*}
In particular, if $(M,\varphi,\xi,\eta,g)$ is a paracontact metric $(-1,\mu)$-space, then the deformed manifold is another paracontact metric $(-1,\mu')$-space with $ \mu'=\frac{\mu-2+2c}{c}$.

Therefore, given a $(-1,2)$-space, a  $\mathcal{D}_c$-homothetic deformation with arbitrary $c\neq0$ will give us another paracontact metric $(-1,2)$-space. Given a paracontact metric $(-1,0)$-space, if we $\mathcal{D}_c$-homothetically deform it with $c=\frac2{2-\mu}\neq 0$ for some $\mu\neq2$,  we will obtain a paracontact metric  $(-1,\mu)$-space with $\mu\neq2$. A sort of converse is also possible: given a $(-1,\mu)$-space  with $\mu\neq2$, a  $\mathcal{D}_c$-homothetic deformation with $c=1-\frac{\mu}{2}\neq 0$ will give us a paracontact metric $(-1,0)$-space.
The case $\mu=0$, $h\neq0$  is also special because the manifold satisfies \eqref{parasasakian-r}
but it is not paraSasakian.


Examples of non-paraSasakian paracontact metric $(-1,2)$-spaces of any possible dimension and constant rank of $h$ were presented in \cite{mio}:

\begin{example}[$(2n+1)$-dimensional paracontact metric $(-1,2)$-space with $\text{rank}(h)=m \in \{1,\ldots,n \}$] \label{ex-mu2-hm-n}

Let $\mathfrak{g}$ be the $(2n+1)$-dimensional Lie algebra with basis $\{\xi,X_1, Y_1, \ldots, X_{n}, Y_{n} \}$ such that the only non-zero Lie brackets are:
\begin{align*}
[\xi,X_i]&=Y_i, \quad i=1,\ldots,m, \\
[X_i,Y_j]&=
\begin{cases}
\delta_{ij}(2\xi +\sqrt2 (1+\delta_{im})Y_m)\\
\qquad+(1-\delta_{ij})\sqrt2 (\delta_{im}Y_j+\delta_{jm}Y_i),  & i,j=1,\ldots,m,\\
\delta_{ij}(2\xi+\sqrt2 Y_i),                                                                          & i,j=m+1, \ldots,n, \\
\sqrt2 Y_i,                                                                                             & i=1, \ldots,m, \; j=m+1, \ldots,n.
\end{cases}
\end{align*}

If we denote by $G$ the Lie group whose Lie algebra is $\mathfrak{g}$, we can define a left-invariant paracontact metric structure on $G$ the following way:
\[
\varphi \xi=0,  \quad \varphi X_i=X_i, \quad  \varphi Y_i=-Y_i,  \quad i=1, \ldots,n,
\]
\[
\eta(\xi)=1,  \quad \eta(X_i)=\eta(Y_i)=0, \quad i=1,\ldots,n.
 \]
The only non-vanishing components of the metric are
\[
g(\xi,\xi)=g(X_i,Y_i)=1, \quad i=1,\ldots,n.
\]
A straightforward computation gives that $h X_i=Y_i$ if $i=1,\ldots,m$, $h X_i=0$ if $i=m+1,\ldots,n$ and $h Y_j=0$ if  $j=1, \ldots,n$, so $h^2=0$ and $\text{rank}(h)=m$. Furthermore, the manifold is a $(-1,2)$-space.
\end{example}

Examples of non-paraSasakian paracontact metric $(-1,0)$-spaces of any possible dimension and $\text{rank}(h)=1$ were also given in \cite{mio}:

\begin{example} [$(2n+1)$-dimensional paracontact metric $(-1,0)$-space with  $\text{rank}(h)=1$] \label{ex-mu0-h1}
Let $\mathfrak{g}$ be the $(2n+1)$-dimensional Lie algebra with basis $\{\xi,X_1,Y_1,\ldots,X_n,Y_n \}$ such that the only non-zero Lie brackets are:
\begin{align*}
[\xi,X_1]&=X_1+Y_1,    & [\xi,Y_1]&=-Y_1,    & [X_1,Y_1]&=2 \xi, \\
[X_i,Y_i]&=2(\xi+Y_i), & [X_1,Y_i]&=X_1+Y_1, & [Y_1,Y_i]&=-Y_1, \quad i=2,\ldots,n.
\end{align*}

If we denote by $G$ the Lie group whose Lie algebra is $\mathfrak{g}$, we can define a left-invariant paracontact metric structure on $G$ the following way:
\[
\varphi \xi=0,  \quad \varphi X_1=X_1, \quad  \varphi Y_1=-Y_1,  \quad\varphi X_i=-X_i, \quad  \varphi Y_i=Y_i, \quad i=2, \ldots,n,
\]
\[
\eta(\xi)=1,  \quad \eta(X_i)=\eta(Y_i)=0, \quad i=1,\ldots,n.
 \]
The only non-vanishing components of the metric are
\[
g(\xi,\xi)=g(X_1,Y_1)=1, \quad g(X_i,Y_i)=-1, \quad i=2,\ldots,n.
\]
A straightforward computation gives that $h X_1=Y_1$, $h Y_1=0$ and $h X_i=hY_i=0$, $i=2, \ldots,n$, so $h^2=0$ and $\text{rank}(h)=1$.

Moreover, by basic paracontact metric properties and Koszul's formula we obtain that
\begin{gather*}
\nabla_\xi X_1=0,  \quad \nabla_\xi Y_1=0, \quad \nabla_\xi X_i=X_i, \quad  \nabla_\xi Y_i=-Y_i, \quad i=2,\ldots,n,\\
\nabla_{X_i} Y_1=\delta_{i1}\xi, \quad \nabla_{X_i} Y_j=\delta_{ij} (\xi+2Y_i), \quad \nabla_{Y_1} X_1=-\xi, \quad \nabla_{Y_i} X_j=-\delta_{ij} \xi, \quad i,j=2,\ldots,n,\\
\nabla_{X_1} X_j=0, \quad \nabla_{Y_1} Y_1=\nabla_{Y_1} Y_j=0, \quad \nabla_{Y_j} Y_1=Y_1, \quad i=2,\ldots,n,
\end{gather*}
and thus
\begin{align*}
R(X_i,\xi)\xi &=-X_i, \quad i=1,\ldots,n,\\
R(Y_i,\xi)\xi &=-Y_i, \quad i=1,\ldots,n,\\
R(X_i,X_j)\xi &=R(X_i,Y_j)\xi=R(Y_i,Y_j)\xi=0, \quad i,j=1,\ldots,n.
\end{align*}
Therefore, the manifold is also a $(-1,0)$-space.
\end{example}

To our knowledge, the previous example is the first paracontact metric $(-1,\mu)$-space with $h^2=0$, $h\neq0$ and  $\mu\neq2$ that was constructed in dimensions greater than $3$. For dimension $3$, Example~4.6 of \cite{CP} was already known.

In dimension $5$, there also exist examples of paracontact metric $(-1,0)$-space with $\text{rank}(h)=2$ and in dimension $7$ of $\text{rank}(h)=2,3$, as shown in \cite{mio}. Higher-dimensional examples of paracontact metric $(-1,0)$-spaces with $\text{rank}(h) \geq2$ were not included, which will be remedied in Example~\ref{ex-mu0-h2+}. We will also see how to construct a $3$-dimensional paracontact metric $(-1,0)$-space  and $(-1,2)$-space where the rank of $h$ is not constant.

\section{New proof of Theorem~\ref{th-h}}

We will now present a revised proof of Theorem~\ref{th-h} that does not use \cite{O} when $h\neq0$ but constructs the basis explicitly.

\begin{proof}
Since $\kappa=-1$, we know from \cite{CKM} that $h^2=0$. If  $h=0$, then $R(X,Y)\xi=-(\eta(Y)X-\eta(X)Y)$, for all $X,Y$ on $M$ and $\xi$ is a Killing vector field, so Theorem~3.1 of \cite{mio} gives us that the manifold is paraSasakian.

If $h\neq0$, then let us take a point $p\in M$ such that $h_p \neq0$. We know that $\xi$ is a global vector field such that $g(\xi,\xi)=1$, that $h \xi=0$ and that $h$ is self-adjoint, so $Ker \eta_p$ is $h$-invariant and $h_p : Ker \eta_p \mapsto Ker \eta_p$ is a non-zero linear map such that $h_p^2=0$.
We will now construct a basis $\{ X_1,Y_1,\ldots,X_n,Y_n \}$ of $Ker \eta_p$ that satisfies all of our requirements.

Take a non-zero vector $v \in Ker \eta_p$ such that $h_p v\neq0$, which we know exists because $h_p\neq0$. Then we write $Ker \eta_p=L_1 \oplus L_1^\perp$, where $L_1=\langle v, h_p v \rangle$. Both $L_1$ and $L_1^\perp$ are $h_p$-invariant because $h_p$ is self-adjoint. Moreover, $g_p (v, h_p v) \neq0$ because $g_p (h_p v, h_p v)=0=g_p (h_p v, w)$ for all $w \in L_1^\perp$, $h_p v \neq0$ and $g$ is a non-degenerate metric.We now distinguish two cases:
\begin{enumerate}
\item If $g_p(v,v)=0$, then we can take $X_i=\frac{1}{\sqrt{|g_p(v,h_p v)|}} v$ and $Y_i=\frac{1}{\sqrt{|g_p(v,h_p v)|}}h_p v$, which  satisfy that $g_p(X_i,X_i)=0=g_p(Y_i,Y_i)$,  $g_p(X_i,Y_i)=\pm 1$ and $h_p X_i=Y_i$.
\item If $g_p(v,v)\neq0$, then $v'=v-\frac{g_p(v,v)}{g_p(v,h_pv)} h_p v$ satisfies that $g_p(v',v')=0$, so we can take $X_i=\frac{1}{\sqrt{|g_p(v',h_p v)'|}} v'$, $Y_i=\frac{1}{\sqrt{|g_p(v',h_p v')|}}h v'$. We have again that $g_p(X_i,X_i)=0=g_p(Y_i,Y_i)$, $g_p(X_i,Y_i)=\pm 1$  and $h_p X_i=Y_i$.
\end{enumerate}

In both cases, $L_1=\langle X_i,Y_i\rangle$, so we now take a non-zero vector  $v \in L_1^\perp$ and check if $h_p v \neq0$. We know that we can take $v$ such that $h_p v\neq0$ in this step as many times as the rank of $h_p$, which is at minimum $1$ (since $h_p\neq0$) and at most $n$ because $\dim Ker \eta_p=2n$ and the spaces $L_1$ have dimension $2$.

If we denote by $m$ the rank of $h_p$, then we can write $Ker \eta_p$ as the following direct sum of mutually orthogonal subspaces:
\[
Ker \eta_p=L_1 \oplus L_2 \oplus \cdots \oplus L_m \oplus V=\langle X_1, Y_1, \ldots, X_m,Y_m \rangle \oplus V,
\]
where $h_p v=0$ for all $v \in V$. Each $L_i$ is of signature $(1,1)$ because $\{ \tilde{X}_i=\frac1{\sqrt2} (X_i+Y_i),\tilde{Y}_i=\frac1{\sqrt2} (X_i-Y_i) \}$ is a pseudo-orthonormal basis such that $g_p(\tilde{X}_i,\tilde{X}_i)=-g_p(\tilde{Y}_i,\tilde{Y}_i)=g_p(X_i,Y_i)=\pm 1$, $g_p(\tilde{X}_i,\tilde{Y}_i)=0$. Therefore, $\langle X_1, Y_1, \ldots, X_m,Y_m \rangle $ is of signature $(m,m)$ and, since $Ker \eta_p$ is of signature $(n,n)$, we can take a pseudo-orthonormal basis
$\{ v_1, \ldots, v_{n-m}, w_1, \ldots, w_{n-m} \}$ of $V$ such that $g_p(v_i,v_j)= \delta_{ij}$ and $g_p(w_i,w_j)= -\delta_{ij}$. Therefore, it suffices to define $X_{m+i}=\frac1{\sqrt2} (v_i+w_i), Y_{m+i}=\frac1{\sqrt2} (v_i-w_i)$ to have $g_p(X_i,X_i)=0=g_p(Y_i,Y_i)$,  $g_p(X_i,Y_i)= 1$ and $h_p X_i=h_p Y_i=0$, $i=m+1, \ldots, n$.

If $n=1$, then $\varphi_p X_1=\pm X_1$ and $\varphi_p Y_1=\mp Y_1$ follow directly from paracontact metric properties and the definition of the basis $\{ X_1,Y_1,\ldots,X_n,Y_n \}$.
\end{proof}

It is worth mentioning that Theorem~\ref{th-h} is true only pointwise, i.e. $\text{rank}(h_p)$ does not need to be the same for every $p\in M$. Indeed, we will see in Examples~\ref{ex-mu2-non-constant} and \ref{ex-mu0-non-constant} that we can construct paracontact metric $(-1,\mu)$-spaces such that $h$ is zero in some points and non-zero in others.

\section{New examples}\label{sec-examples}

We will first present an example of $(2n+1)$-dimensional paracontact metric $(-1,0)$-space with  rank of $h$ greater than $1$. This means that, together with Examples~\ref{ex-mu2-hm-n} and \ref{ex-mu0-h1}, we have examples of paracontact metric $(-1,\mu)$-spaces of every possible dimension and constant rank of $h$ when $\mu=0$ and $\mu=2$.

\begin{example} [$(2n+1)$-dimensional paracontact metric $(-1,0)$-space with  $\text{rank}(h)=m\in \{ 2, \ldots, n\}$] \label{ex-mu0-h2+}
Let $\mathfrak{g}$ be the $(2n+1)$-dimensional Lie algebra with basis $\{\xi,X_1,Y_1,\ldots,X_n,Y_n \}$ such that the only non-zero Lie brackets are:
\begin{align*}
[\xi,X_1]&=X_1+X_2+Y_1,                      &  [\xi,Y_1]&=-Y_1+Y_2,\\
[\xi,X_2]&=X_1+X_2+Y_2,                      &  [\xi,Y_2]&=Y_1-Y_2,\\
[\xi,X_i]&=X_i+Y_i,  \quad i=3,\ldots,m,     &  [\xi,Y_i]&=-Y_i, \quad i=3,\ldots,m,
\end{align*}
\begin{align*}
[X_i,X_j]&=
\begin{cases}
\sqrt{2} X_1,       &\text{ if } i=1, \; j=2,\\
-\sqrt2 X_j         &\text{ if } i=2, \; j=3\ldots,m,\\
\sqrt2 [\xi,X_i],   &\text{ if } i=1,\ldots,m, \; j=m+1,\ldots,n,
\end{cases}
\\
[Y_i,Y_j]&=
\begin{cases}
\sqrt{2} (-Y_1+Y_2),          &\text{ if } i=1,   \; j=2,\\
\sqrt2 Y_j,                   &\text{ if } i=1,2, \; j=3,\ldots,m,
\end{cases}
\\
[X_i,Y_i]&=
\begin{cases}
2\xi+\sqrt2 (X_2+Y_2)          &\text{ if } i=1,\\
-2\xi+\sqrt2 X_1,              &\text{ if } i=2,\\
-2\xi+\sqrt2 (X_1-X_2-Y_2),    &\text{ if } i=3,\ldots,m,\\
-2\xi-\sqrt2 X_i,              &\text{ if } i=m+1,\ldots,n,
\end{cases}
&
\\
\underset{\mbox{$i\neq j$}}{[X_i,Y_j]}&=
\begin{cases}
\sqrt2 (Y_1+X_2)         &\text{ if } i=1,   \; j=2,\\
\sqrt2 X_1,              &\text{ if } i=2,   \; j=1,\\
\sqrt2 X_j,              &\text{ if } i=1,2, \; j=3,\ldots,m,\\
\sqrt2 Y_i,              &\text{ if } i=3,\ldots,m,   \; j=2,\\
-\sqrt2 [\xi,Y_j],       &\text{ if } i=m+1,\ldots,n, \; j=1,\ldots,m.
\end{cases}
\end{align*}

If we denote by $G$ the Lie group whose Lie algebra is $\mathfrak{g}$, we can define a left-invariant paracontact metric structure on $G$ the following way:
\[
\varphi \xi=0,  \quad \varphi X_i=X_i, \quad  \varphi Y_i=-Y_i,  \quad i=1, \ldots,n,
\]
\[
\eta(\xi)=1,  \quad \eta(X_i)=\eta(Y_i)=0, \quad i=1,\ldots,n.
 \]
The only non-vanishing components of the metric are
\[
g(\xi,\xi)=g(X_1,Y_1)=1, \quad g(X_i,Y_i)=-1, \quad i=2,\ldots,n.
\]
A straightforward computation gives that $h X_i=Y_i$, $i=1, \ldots,m$, $h X_i=0$, $i=m+1, \ldots,n$ and  $h Y_i=0$,  $i=1, \ldots,n$, so $h^2=0$ and $\text{rank}(h)=m$.

Moreover, very long but direct computations give that
\begin{align*}
R(X_i,\xi)\xi &=-X_i, \quad i=1,\ldots,n,\\
R(Y_i,\xi)\xi &=-Y_i, \quad i=1,\ldots,n,\\
R(X_i,X_j)\xi &=R(X_i,Y_j)\xi=R(Y_i,Y_j)\xi=0, \quad i,j=1,\ldots,n.
\end{align*}
Therefore, the manifold is also a $(-1,0)$-space.
\end{example}

\begin{remark}
Note that the previous example is only possible when $n\geq2$. If $n=1$, then we can only construct examples of $\text{rank}(h)=1$, as in Example~\ref{ex-mu0-h1}.

In the definition of the Lie algebra of the previous example, some values of $i$ and $j$ are not possible for $m=2$ or $m=n$. In that case, removing the affected Lie brackets from the definition will give us valid examples nonetheless.
\end{remark}

We will present now an example of $3$-dimensional paracontact metric $(-1,2)$-space and one of $3$-dimensional paracontact metric $(-1,0)$-space, such that $\rank(h_p)=0$ or $1$ depending on the point $p$ of the manifold. These are the first examples of paracontact metric $(\kappa,\mu)$-spaces with $h$ of non-constant rank that are known.

\begin{example}[$3$-dimensional paracontact metric $(-1,2)$-space with $\text{rank}(h_p)$ not constant] \label{ex-mu2-non-constant}

We consider the manifold $M=\mathbb{R}^3$ with the usual cartesian coordinates $(x,y,z)$. The vector fields
\[
e_1= \frac{\partial}{\partial x}+x z \frac{\partial}{\partial y}-2y\frac{\partial}{\partial z}, \quad
 e_2=\frac{\partial}{\partial y}, \quad
 \xi=\frac{\partial}{\partial z}
\]
are linearly independent at each point of $M$. We can compute
\[
[e_1,e_2]=2 \, \xi, \quad
[e_1,\xi]=-x \, e_2, \quad
[e_2,\xi]=0.
\]

We define the semi-Riemannian metric $g$ as the non-degenerate one whose only non-vanishing components are $g(e_1,e_2)=g(\xi,\xi)=1$, and the $1$-form $\eta$ as $\eta=2y dx+dz$, which satisfies $\eta(e_1)=\eta(e_2)=0$, $\eta(\xi)=1$. Let $\varphi$ be the $(1,1)$-tensor field defined by $\varphi e_1=e_1, \varphi e_2=-e_2, \varphi \xi=0$. Then
\begin{align*}
d\eta(e_1,e_2)&=\frac12 (e_1(\eta(e_2))-e_2(\eta(e_1))- \eta([e_1,e_2]))=-1=-g(e_1,e_2)=g(e_1,\varphi e_2),\\
d\eta(e_1,\xi)&=\frac12 (e_1(\eta(\xi))-\xi(\eta(e_1))- \eta([e_1,\xi])=0=g(e_1,\varphi \xi),\\
d\eta(e_2,\xi)&=\frac12 (e_2(\eta(\xi))-\xi(\eta(e_2))- \eta([e_2,\xi])=0=g(e_2,\varphi \xi).
\end{align*}
Therefore, $(\varphi,\xi,\eta,g)$ is a paracontact metric structure on $M$.

Moreover, $h\xi=0$, $h e_1=x  e_2$, $h e_2=0$. Hence, $h^2=0$ and, given $p=(x,y,z) \in \mathbb{R}^3$, $\rank(h_p)=0$ if $x=0$ and $\rank(h_p)=1$ if $x\neq0$.

Let $\nabla$ be the Levi-Civita connection. Using the properties of a paracontact metric structure and Koszul's formula
\begin{equation}\label{eq-koszul}
2g(\nabla_X Y,Z) = X(g(Y, Z)) + Y(g(Z, X)) - Z(g(X, Y)) - g(X, [Y, Z]) - g(Y, [X, Z]) + g(Z, [X, Y]),
\end{equation}
we can compute
\[
\nabla_\xi \xi=0, \quad \nabla_{e_1} \xi=-e_1-x e_2, \quad \nabla_{e_2} \xi=e_2, \quad
\nabla_\xi e_1=-e_1, \quad \nabla_\xi e_2=e_2,
\]
\[
\nabla_{e_1} e_1=x \xi, \quad
\nabla_{e_2} e_2=0, \quad\nabla_{e_1} e_2= \xi, \quad \nabla_{e_2} e_1 =-\xi.
\]
Using the following definition of the Riemannian curvature
\begin{equation}\label{eq-Riemann}
R(X,Y)Z=\nabla_X \nabla_Y Z-\nabla_Y \nabla_X Z-\nabla_{[X,Y]} Z,
\end{equation}
we obtain
\[
R(e_1,\xi)\xi=-e_1+2 h e_1, \quad  R(e_2,\xi)\xi=-e_2+2 h e_2, \quad R(e_1,e_2)\xi=0,
\]
so the paracontact metric manifold $M$ is also a $(-1,2)$-space.
\end{example}

\begin{remark}
The previous example does not contradict Theorem~\ref{th-h}, as we will see by constructing explicitly the basis of the theorem on each point $p$ where $h_p\neq0$, i.e., on every point $p=(x,y,z)$ such that $x\neq0$.

Indeed, let us take a point $p=(x,y,z) \in \mathbb{R}^3$. If $x\neq0$, then we define $X_1=\frac{{e_1}_p}{\sqrt{|x|}}$, $Y_1=\frac{h_p{e_1}_p}{\sqrt{|x|}}$. We obtain that $\{ \xi_p, X_1, Y_1 \}$ is a basis of $T_p(\mathbb{R}^3)$ that satisfies that:
\begin{itemize}
\item the only non-vanishing components of $g$ are
$
g_p(\xi_p,\xi_p)=1, \quad g_p(X_1,Y_1)=sign(x),
$

\item
the tensor $h$ can be written as
$
{h_p}_{| \langle\xi_p, X_1,Y_1 \rangle}=
\begin{pmatrix}
0 & 0 & 0\\
0 & 0 & 0\\
0 & 1 & 0
\end{pmatrix},
$

\item $\varphi_p \xi=0, \quad \varphi_p X_1 =  X_1, \quad \varphi_p Y_1 =-  Y_1.$
\end{itemize}

\end{remark}

\begin{example}[$3$-dimensional paracontact metric $(-1,0)$-space with $\text{rank}(h_p)$ not constant] \label{ex-mu0-non-constant}

We consider the manifold $M=\mathbb{R}^3$ with the usual cartesian coordinates $(x,y,z)$. The vector fields
\[
e_1= \frac{\partial}{\partial x}+xe^{-2z} \frac{\partial}{\partial y}
-2y\frac{\partial}{\partial z}, \quad
 e_2=\frac{\partial}{\partial y}, \quad
 \xi=\frac{\partial}{\partial z}
\]
are linearly independent at each point of $M$. We can compute
\[
[e_1,e_2]= 2 \, \xi, \quad
[e_1,\xi]=2x e^{-2z}  \, e_2, \quad
[e_2,\xi]=0.
\]

We define the semi-Riemannian metric $g$ as the non-degenerate one whose only non-vanishing components are $g(e_1,e_2)=g(\xi,\xi)=1$, and the $1$-form $\eta$ as $\eta=2y dx+dz$, which satisfies $\eta(e_1)=\eta(e_2)=0$, $\eta(\xi)=1$. Let $\varphi$ be the $(1,1)$-tensor field defined by $\varphi e_1=e_1, \varphi e_2=-e_2, \varphi \xi=0$. Then
\begin{align*}
d\eta(e_1,e_2)&=\frac12 (e_1(\eta(e_2))-e_2(\eta(e_1))- \eta([e_1,e_2]))=-1=-g(e_1,e_2)=g(e_1,\varphi e_2),\\
d\eta(e_1,\xi)&=\frac12 (e_1(\eta(\xi))-\xi(\eta(e_1))- \eta([e_1,\xi])=0=g(e_1,\varphi \xi),\\
d\eta(e_2,\xi)&=\frac12 (e_2(\eta(\xi))-\xi(\eta(e_2))- \eta([e_2,\xi])=0=g(e_2,\varphi \xi).
\end{align*}
Therefore, $(\varphi,\xi,\eta,g)$ is a paracontact metric structure on $M$.

Moreover, $h\xi=0$, $h e_1=-2x e^{-2z} e_2$, $h e_2=0$. Hence, $h^2=0$ and, given $p=(x,y,z) \in \mathbb{R}^3$, $\rank(h_p)=0$ if $x=0$ and $\rank(h_p)=1$ if $x\neq0$.

Let $\nabla$ be the Levi-Civita connection. Using the properties of a paracontact metric structure and Koszul's formula \eqref{eq-koszul}, we can compute
\[
\nabla_\xi \xi=0, \quad \nabla_{e_1} \xi=-e_1+2x e^{-2z} e_2, \quad \nabla_{e_2} \xi=e_2, \quad
\nabla_\xi e_1=-e_1, \quad \nabla_\xi e_2=e_2,
\]
\[
\nabla_{e_1} e_1=-2xe^{-2z} \xi, \quad
\nabla_{e_2} e_2=0, \quad\nabla_{e_1} e_2=  \xi, \quad \nabla_{e_2} e_1 =-\xi.
\]
Using now \eqref{eq-Riemann}, we obtain
\[
R(e_1,\xi)\xi=-e_1, \quad  R(e_2,\xi)\xi=-e_2, \quad R(e_1,e_2)\xi=0,
\]
so the paracontact metric manifold $M$ is also a $(-1,0)$-space.
\end{example}

\textbf{Acknowledgements.}
The author would like to thank Prof. Mart\'{\i}n Avenda\~no for his invaluable help.

\end{document}